\documentclass[a4paper,11pt]{article}
\usepackage{pstricks}
\usepackage{amsmath}
\usepackage{amsfonts}
\usepackage{pstcol}
\usepackage{amscd}
\usepackage[all]{xy}

\font\msbm=msbm10
\newcommand{\B}[1]{\hbox{\msbm #1}}

\newcommand{\al}{{\alpha}}

\newcommand{\be}{{\beta}}

\newcommand{\Om}{{\Omega}}
\newcommand{\om}{{\omega}}

\newcommand{\De}{{\Delta}}

\newcommand{\la}{{\lambda}}
\newcommand{\La}{{\Lambda}}
\newcommand{\si}{{\sigma}}
\newcommand{\Si}{{\Sigma}}

\newcommand{\qed}{\hfill {\bf QED}\medskip}

\newcommand{\map}[1]{\stackrel {#1}\longrightarrow}

\newcommand{\Mo}{(M,\omega )}

\newcommand{\hook}{\hookrightarrow}
\newcommand{\cp}{\B {CP}}
\newcommand{\pf}{\NI {\bf Proof: }}

\newcommand{\BS}{{\bigskip}}
\newcommand{\NI}{{\noindent}}

\newcommand{\QED}{\hfill {\bf QED}\medskip}

\newtheorem{theorem}{Theorem}[subsection]
\newtheorem{thm}[theorem]{Theorem}
\newtheorem{cor}[theorem]{Corollary}
\newtheorem{defin}[theorem]{Definition}
\newtheorem{rem}[theorem]{Remark}
\newtheorem{lemma}[theorem]{Lemma}
\newtheorem{prop}[theorem]{Proposition}
\newtheorem{ex}[theorem]{Example}

\numberwithin{equation}{section}

\title{Characteristic classes of smooth fibrations}
\author{Tadeusz Januszkiewicz
        \thanks {Partially supported by
        State Committee for Scientific Research grant 2 P03A 035 20} \\
         Wroc\l aw University and IMPAN
        \and Jaros\l aw K\c edra
        \thanks{Partially supported by State Committee for Scientific
        Research grant 5 P03A 017 20 }\\
        Max Planck Institut,
        Szczecin University and \\
        IM PAN}

\date{\today}

\begin{document}

\psset{unit=1.5cm}

\maketitle

\begin{abstract}
We construct characteristic classes of
smooth (Hamiltonian) fibrations as as fiber
integrals of products of Pontriagin (or Chern) classes of vertical
vector bundles over the total space of the universal fibration.
We give explicit formulae of these fiber integrals for toric manifolds
and get estimates of the dimension of the
cohomology groups of classifying spaces.
{\bf AMS classification(2000):} Primary 55R40; Secondary 57R17
\end{abstract}

\tableofcontents

\BS
\NI
{\bf Acknowledgments:} We would like to thank Dusa McDuff for useful comments
on the earlier version of this paper and Friedrich Hirzebruch for calling our
attention to multiplicativity properties of genera.

\BS

\section{Background}\label{S:back}

\BS

\subsection{Universal fibrations and characteristic classes}
\label{SS:univ}

Let $G$ be a topological group for which the universal principal
fibration $G\to E_G\to B_G$ is locally trivial.
Suppose $G$ acts (smoothly) on a manifold $M$.
Consider the associated bundle
$$M\hook M_G = E_G\times _G M \to B_G.$$

\NI
Any fibration       $M\hook P\to B$            with structural
group $G$ and fiber $M$ is obtained as a pullback of $M_G$ for an appropriate
map $f:B\to BG$

$$
\begin{CD}\label{D:pull}
             P=f^*M_G    @>{\bar f}>>   M_G \\
             @V{\pi _B}VV              @V{\pi }VV   \\
             B           @>{f}>>         BG.
\end{CD}
$$

\NI
A {\bf  $\chi $-characteristic class}
of  the $G$-fibration  $M\hook P\to B$
is a cohomology class of the form
$f^*\chi \in H^*(B)$, where $\chi \in H^*(BG)$.

\subsection{Equivariant  bundles}\label{SS:vert}

A $G$-equivariant bundle $\xi \to M$
extends to a  bundle $\xi_G\to M_G$ such
that the following diagram commutes

$$
\begin{CD}
\xi    @>>>  \xi _G    @>>>    B_G          \\
 @VVV          @VVV           @V\hbox{Id}VV  \\
 M     @>i>>   M_G      @>>>    B_G.
\end{CD}
$$

\noindent
More precisely, $\xi _G= (\xi \times E_G) \slash G$ and $i^*\xi _G = \xi $, where
$i:M\to M_G$  is an inclusion of the fiber.
A particularly simple case is when $\xi = TM$ is the tangent bundle of $M$,
and the $G$ action is smooth.
Then we call $TM_G$ the {\bf vertical bundle}
of the fibration $M_G\to B_G$.

We will need additional structures in these extended
bundles, the most important being complex structure.
If $G$ acts on a vector bundle preserving a complex
structure $J$, then $\xi _G$ admits a complex structure.
We can relax somewhat invariance condition of
$J$ with respect to $G$ and
still obtain a complex structure in $\xi_G$.

\BS
\begin{prop} \label{P:cs}
Let $\xi \to M$ be a $G$-equivariant real vector bundle.
Suppose that there exists a contractible $G$-invariant set $\cal J$
of complex structures in $\xi $. Then there exists a complex
structure $J_G$  in $\xi _G$
whose restriction to each fiber is homotopic to any $J\in \cal J$.
In particular,
$i^*c_k(\xi _G,J_G) = c_k(\xi ,J)$, for $J\in \cal J$.
\end{prop}

\BS
\NI
{\bf Proof:}
Consider the map $\xi\times E_G\times {\cal J}
\to M\times E_G\times {\cal J}$.
This is clearly equivariant with respect to the $G$ action,
and thus defines a projection of a bundle

$$\xi^J_G = \xi\times E_G\times {\cal J}\slash G
\to M\times E_G\times {\cal J}\slash G.$$

\noindent
Since $\cal J$ is contractible, the base is homotopy equivalent to
$M_G$ and  $\xi^J_G$ is essentially $\xi_G$, that is it is a pullback by the
inverse of the
homotopy equivalence $E_G\times {\cal J}\slash G \to E_G\slash G$.
The bundle $\xi\times {\cal J}\to M\times {\cal J}$
admits a tautological complex structure $\widetilde J$:

$$\widetilde J_{(m,J)}(v)=J(m)(v),$$

\noindent
which is $G$-equivariant and thus define a complex
structure in $\xi_G^J$.
Pulling it back we obtain a required complex structure in $\xi_G$.

\QED

Proposition \ref{P:cs} is of very general nature
and should have many variants for groups acting on vector bundles
and preserving geometric structures.
The motivating example for the Proposition \ref{P:cs} comes
from symplectic geometry. It is a remark of fundamental importance,
due to Gromov,
that the set of complex structures $J$ tamed by a symplectic form
$\om $ (i.e. such that $\om (J\circ , \circ )$ is positive
 definite symmetric form) is contractible.
It
consists of sections of the bundle over $M$, which is
associated to the frame bundle and has a contractible fiber
$Sp(2n;\B R)\slash U(n)$.

\section{Characteristic classes from fiber integrals}
\label{SS:fi}

Recall that for any fibration $M\hook P\to B$, where $M$ is
oriented $m$-dimensional compact manifold,
there is a homomorphism of $H^*(B)$-modules

$$\pi _*:H^*(P)\to H^{*-m}(B)$$

\NI
called {\bf fiber integration}
\footnote{In case when $M$ is not compact $\pi _*$ is defined
on cohomology with compact support.}
(see [AB],[GS Chapter 10] for detailed description).
Although fiber integration behaves badly with respect to
the cup product, it has useful properties.
The most important to us is its naturality.

Let $M\hook P\map{\pi } B$ be a fibration with $m$-dimensional
fiber and let be given its pullback

$$
\CD
        f^*P  @>{\tilde f}>>  P \\
        @VV\pi ^{\prime }V  @VV\pi V \\
        B^{\prime } @>f>>  B.
\endCD
$$

\NI
Then the naturality of fiber integration means that the following
diagram commutes

$$
\CD \label{D:nat}
      H^*(f^*P) @<\tilde f^*<<  H^*(P) \\
       @VV\pi ^{\prime }_*V  @VV\pi _*V \\
      H^{*-m}(B^{\prime }) @<f^*<< H^{*-m}(B).
\endCD
$$

\BS

We use fiber integration 
to construct characteristic classes of $G$-bundles
from characteristic classes of $\xi _G$-bundles.
This construction has several variants depending on
the structure of $\xi _G$ which we describe below. 
The simple characteristic classes on $M_G$
become fairly involved when pushed down. 
This of course agrees with the 
general philosophy of reduction.

\subsection{The Pontriagin classes}\label{E:pont}

Let $G=Dif\! f(M)$ and
suppose $\xi $ is a 
$G$-equivariant vector bundle. 
For example, $\xi $ might be the tangent bundle.
Take $\xi _G\to M_G$
and a monomial in Pontriagin classes 
$p_I(\xi _G):=p_{i_1}(\xi _G)\cup \dots \cup p_{i_k}(\xi _G)$. 
Then integrate it over the fiber of the fibration
$M\to M_G\to B_G$, to obtain classes 
$\pi _*^G(p_I(\xi _G))\in H^*(B_G)$.
This is depicted in the following diagram where
$f_{\xi _G}:M_G \to B_{O(n)}$ is a classifying
map and the curved arrow is the composition.

$$
\xymatrix
{     H^*(M_G) \ar[d]^{\pi _*} &
      H^*(B_{O(n)})  \ar[l]_{f_{\xi _G}} \ar@/_1.1pc/[dl] \\
      H^{*-n}(B_G) 
}
$$

\NI

\subsection{The Pontriagin and the Euler classes}\label{E:euler}

If $G\subset Dif\! f(M)$ 
acts on an oriented bundle $\xi $ preserving the orientation 
then we can also consider the Euler class of $\xi $. Integrating 
monomials in  the Euler and Pontriagin classes along the fiber we get elements in 
$H^*(B_G)$. Again, if $\xi $ is the tangent bundle of an oriented
manifold $M$ then we obtain an information about
the cohomology of the classifying space of the group of orientation
preserving diffeomorphisms.

\BS
For an orientable 2-dimensional manifold $M$ the classes
$\pi _*(eu(TM_G)^k)$ are the, so called, Miller-Morita-Mumford
classes and have been extensively studied \cite{mo}.

\subsection{The Chern classes}\label{E:chern}

Let $G\subset Dif\! f(M)$.
Suppose now that $\xi $ is an $G$-equivariant
vector bundle over $M$ admitting a family of complex structures
as in Proposition \ref{P:cs}.
Consider a monomial in the Chern classes
$c_I:=c_{i_1}(\xi _G)\cup \dots \cup c_{i_k}(\xi _G)$, and
obtain classes
$\pi _*^G(c_I)\in H^*(B_G)$. If $\xi $ is a tangent bundle
of a symplectic manifold $\Mo $, $G=Symp\Mo $ 
and $\cal J$ is the set of tamed almost complex structures,
then we get characteristic classes 
of symplectic fibrations.

\subsection{Hamiltonian actions and coupling class}\label{SS:ha}

It is also  interesting to consider the group
$Ham\Mo \subset Symp\Mo $ of Hamiltonian symplectomorphisms.
Then  an additional characteristic
class, arises [GS, Section 9.5],[Po, page 70]:

\begin{defin}\label{D:cc}
An element
$\Om \in H^2(M_{Ham\Mo })$
satisfying the 
following conditions:
       \begin{enumerate}
             \item
             $\pi _*(\Om ^{n+1}) = 0 $
             \item
              $i^*(\Om ) = [\om ].$
       \end{enumerate}
\noindent
is called the {\bf coupling class}
\end{defin}

\begin{prop}\label{P:cc}
The coupling class is unique and it is equal to

\begin{equation}\label{E:coupling}
\Om := \widetilde{\om } -{1\over (n+1)Vol\Mo }
\pi^*(\pi _*\widetilde{\om }^{n+1}),
\end{equation}

\noindent
where
$\widetilde{\om } \in H^*(M_{Ham\Mo })$ restricts to the
class of symplectic form of the fiber (i.e. 
$i^*(\widetilde{\om }) = [\om ]\in H^*(M)$) and
$Vol\Mo = {1\over n!}\int_M\om ^n$ is the symplectic volume
of $\Mo $.
\end{prop}

\BS
\NI
{\bf Proof:}
Since the universal fibration $M_{Ham\Mo }\to BHam\Mo $
is Hamiltonian, then there is a class
$\widetilde {\om }\in H^2(M_{Ham\Mo })$
which restricts to the class of
symplectic form of the fiber \cite{ms}. 
The class given by the above formula
satisfies the conditions $(1)$ of Definition \ref{D:cc}. 

Since
$H^2(M_{Ham\Mo })\cong H^2(B_{Ham\Mo })\oplus H^2(M)$ \cite{lm},
then if there were another class $\Om _1$ satisfying
the definition then we would have $\Om _1 = \Om + \pi ^*\al$,
for some $\al \in H^2(B_{Ham\Mo })$ and

\begin{eqnarray*}
\pi _*(\Om _1^{n+1}) &=&\pi _*(\Om + \pi ^*\al)^{n+1} \\
& = &\pi _*(\Om ^{n+1} + (n+1)\Om ^n\pi^*\al +...)\\
& = &(n+1)\al \neq 0.
\end{eqnarray*}

\QED

\BS

Let $G= Ham\Mo $ be a the group of Hamiltonian 
symplectomorphisms.
Then integrating monomials in the Chern classes
and the coupling class we obtain an information about
$H^*(B_{Ham\Mo })$ in the same way as in the previous
constructions. It is worth noting that taking the coupling class
into account is essential as we shall see in  Section \ref{SS:ex}.

Reznikov \cite{re} has constructed characteristic classes
of Hamiltonian fibration using the Chern-Weil theory.
As remarked by McDuff
in \cite{m} the classes obtained by Reznikov are fiber integrals of
powers of the coupling class.

\BS

\section{The setup for the computations}

\BS
We are interested in nonvanishing and independence of classes
in $H^*(B_{Dif\! f(M)})$ ($H^*(B_{Symp\Mo })$ or
$H^*(B_{Ham\Mo })$). We study this question
by pulling these classes back to $H^*(B_H)$ using
smooth (symplectic or Hamiltonian) actions of a compact Lie group $H$.
Notice that it suffices to restrict to the case of torus actions,
as $H^*(B_H)=H^*(B_T)^W\subset H^*(B_T)$, where $T\subset H$ is a
maximal torus and $W$ is the Weyl group.

\subsection{The detection function}\label{SS:ex}

In the sequel the group $G$ will denote the group of Hamiltonian
symplectomorphisms or the group of all diffeomorphisms of 
a symplectic manifold $\Mo $. 
It is worth of noting that we obtain more information 
if we take into account
several torus actions. This is the reason for the
following definition.

\begin{defin}\label{D:d}
Let $f_i:T\to G\subset Dif\!f(M) $, $i=1,...,m$ 
be torus actions. 
We define the {\bf detection function} by  

$$
\CD
F:H^*(B_G) \to  \bigoplus _{i=1}^m H^*(B_T)\\
F(\al ):= [f_1(\al ),...,f_m(\al )]
\endCD
$$
\end{defin}

\BS

\NI
We compute the image of the detection function 
on the fiber integrals of characteristic
classes which allows to estimate the Betti numbers of $B_{Ham\Mo }$
and $B_{Dif\!f(M)}$.

\subsection{Localization formula}\label{SS:loc}

Fiber integration for torus actions can be
computed using the Atiyah-Bott-Berline-Vergne \cite{ab,bv}

$$\pi_*(\al ) = \sum_P \pi ^P_*\left( {i^*_P\al \over E(\nu _P)}\right),$$

\noindent
for $\al \in H^*(M_{T})$.
Here $P$ is a connected component of the set of fixed points,
$i_P:P\hook  M$ is the inclusion,
$\pi_*^P:H^*(P_{T})\cong H^*(B_{T})\otimes  H^*(P)\to H^*(B_{T})$ is
fiber integration for the trivial action of $T$ on $P$ and
$E(\nu _P)$  is the equivariant Euler  class of the normal bundle
$\nu _P$ to $P$ in $M$.

In fact the above formula describes the fiber integration
between the localized rings (where the localization is
in the ideal generated by the Euler classes of the normal bundles
of the connected components of the fixed point set). We apply
this to the elements of $H^*(M_{T})$ due to the fact that
the fiber integration commutes with the localization as
depicted in the following diagram, where the right hand
side column consists of localized rings:

$$
\begin{CD}
H^*(M_{T})     @>>>      H^*(M_{T})_{\bullet }\\
@V\pi _*VV                 @V(\pi _*)_{\bullet }VV\\
H^*(B_{T})     @>>>      H^*(B_{T})_{\bullet }
\end{CD}
$$

\subsection{Symplectic toric varieties, moment maps and Delzant polytopes}
\label{S:tor}

\BS

\begin{defin}
{\bf A symplectic toric variety} is a symplectic manifold 
$(M^{2n},\om ) $
endowed with an effective Hamiltonian action of an n-dimensional
torus $T$. 
\end{defin}

Let $a:G\to Ham\Mo $ be a Hamiltonian action of a Lie group $G$ 
on a symplectic manifold $M$. Given $X\in Lie(G)$ a fundamental
vector field $\overline X$ on $M$ is Hamiltonian, which means that

$$\iota _{\overline X}\om = -dH_X,$$

\noindent
for some function $H_X : M\to \B R$. 
This function, initially determined up to a constant,
is chosen so that the map $Lie(G) \to C^{\infty }(M)$
given by $X\to H_X$ is a homomorphism of Lie algebras. The Lie
algebra structure on functions is given by the Poisson bracket.

Then one defines 
{\bf a moment map} $\Phi :M \to  Lie(G)^*$ by

$$\Phi(p)(X) = H_X(p).$$

\BS
\NI
It is a fundamental result of Atiyah,
that for actions of tori the image of a moment map is a convex
polytope $\Delta \subset (\B R^n)^* \cong Lie(T^n)^*$.
Moreover $\Delta$ satisfies the following conditions:
\begin{enumerate}
       \item There  are $n$ edges meeting at each vertex $p\in \Delta $
             (this means that $\Delta $ is a simple polytope).

       \item Every edge including $p$ is of the form $p + tv_i$,
             where $v_i\in (\B Z^n)^*$.

       \item $v_1,...,v_n$ in (2) can be chosen to be a basis of $(\B Z^n)^*$.
\end{enumerate}

\noindent
Such polytope is called {\bf Delzant}. 
Given a Delzant polytope $\Delta $ there is a symplectic 
manifold equipped with a Hamiltonian torus action such that the 
image of the moment map is exactly $\Delta $ [De],[G, Theorem 1.8].
It is useful to describe $\De $ by a system of inequalities of
the form

\begin{equation}\label{E:dp}
\left <x,u_i \right >\geq \la _i,\qquad i-1,...,k,
\end{equation}

\noindent
where $u_i\in \B Z^n$. The vectors $u_i$ can be normalized
by requiring them to be primitive. This normalization together
with inequalities (\ref{E:dp}) determine $u_i$'s uniquely.

Notice that the vectors $u_i$ may be thought as the inward
pointing vectors normal to the faces of the Delzant polytope.
The function which associate the vector $u_i$ to the face
$F_i$ was defined in \cite{dj}(page 423) 
and called the characteristic function.
\BS

\subsection{Equivariant cohomology of toric 
manifolds and the face ring }\label{SS:st}

Consider the universal fibration associated to a toric
manifold $M\hook M_{T} \map{\pi } B_{T}$. We are interested in
the description of the homomorphism 
$\pi ^* :H^*(B_{T}) \to H^*(M_{T})$ induced by the projection in
the universal fibration. By a result  \cite{dj}, the cohomology
ring $H^*(M_{T})$ is isomorphic to the {\bf face ring} of the Delzant
polytope
  \footnote{Face rings are also sometimes called Stanley-Reisner or
            Stanley rings.}.
The face ring $St(\De )$ of a given $n$-dimensional simple
polytope $\De $ is defined to be a graded ring generated by
the $(n-1)$-dimensional faces $F_1,...,F_k$ of $\De $, all generators being 
of degree two, subject to the relations $F_iF_j=0$ iff 
$F_i\cap F_j = \emptyset $ [S, Chapter II].

\begin{rem}
In the sense of the definition given in \cite {s}, our face ring is
the face ring of the simplicial complex dual to $\De $.
\end{rem}

Now we can identify the map $\pi ^*:H^*(B_T)\to H^*(M_T)$.
Since $H^*(B_{T}) \cong \B R[T_1,...,T_n]$ and $H^*(M_{T})\cong St(\De )$
are generated by the elements of degree 2, then it suffices to
describe $\pi ^*$ on the second cohomology. With respect
to the bases $\{T_1,...,T_n\}$ of $H^2(B_{T})$ and 
$\{F_1,...,F_k\}$ of $H^2(M_{T})$, $\pi ^*$ is given by the
matrix $[u_{ij}]$, where $[u_{i1},...,u_{in}] = u_i\in \B Z^n$ 
are the vectors normal to the faces of $\De $ \cite {dj}. 

\BS

Next we describe the Chern classes of the vertical fibration
$TM_{T}$ (see Proposition \ref{P:cs}) in terms of the face
ring $St(\De )$  \cite {dj}.

\begin{prop}\label{P:chern}
The Chern classes $c_i(TM_{T})\in St(\De ) \cong H^*(M_{T})$ 
of the vertical fibration  are given by the following formula:

$$c_i(TM_{T}) = \sum _{1<j_1<...<j_i<k}F_{j_1} ...  F_{j_i}.$$

\end{prop}

{\bf Proof:} 
The idea of the proof is to show that the vertical bundle
$TM_{T}$ is stably isomorphic as a complex bundle to the sum of line bundles,
whose Chern classes are represented by the faces of
the Delzant polytope. This was done by in \cite{dj} in the real
case. The specific choice of orientations of the line bundles
provides the isomorphism preserving complex structures.

To be more precise, we associate a line bundle $L_i \to M_{T}$
to every face $F_i$ of the Delzant polytope in the following way.
The face $F_i$ gives a vector $[u_{i1},...,u_{in}] \in \Bbb Z^n$
as described by the inequalities (\ref{E:dp}). Since 
$\Bbb Z^n \cong \widehat {T^n}$ then we obtain a character
$\chi _i:T^n \to \Bbb S^1$. Now the line bundle $L_i$
is defined as

$$L_i := (M\times E_{T}\times \Bbb C)\slash T,$$

\noindent
where the action of the torus on $\Bbb C$ is given by the
character $\chi _i$. The choice of the {\bf inward pointing} normal 
vectors to the faces of 
the Delzant polytope ensures that the sum of the line bundles
$L_i$ is stably isomorphic to $TM_{T}$ as a complex bundle. 

\QED

Notice that, since the action is Hamiltonian, then we may ask
also about the form of the coupling class in the face ring.

\begin{prop}\label{P:coupling}
The coupling class $\Om \in St(\De )$ is given by 

$$\Om = - \sum \la _i F_i 
        - {1\over (n+1)Vol\Mo }\pi ^*(\pi _*(-\sum \la _i F_i)^{n+1}),$$
where $\la _i\in \B R$ are as in (\ref{E:dp}).
\end{prop}

\noindent
{\bf Proof:} As we have observed in Section \ref{SS:ha}, the coupling
class is given by the formula (\ref{E:coupling}). Thus we have to show
that the class $-\sum \la _i F_i$ restricts to the class of the symplectic
form, that is to say $i^*(-\sum \la _i F_i) = [\om ]$. Here 
$i:M\to M_{T}$ is an embedding of the fiber.
It is known that the class of the symplectic form is equal to
$-\sum \la_1 p_1$, where $p_i\in H^2(M)$ is Poincare dual to
the homology classes given by the preimages of the $(n-1)$-dimensional 
face $F_i$ of the Delzant polytope under the moment map \cite{g}(Appendix 2).
On the other hand, these classes are exactly equal to
$i^*(F_i) $, which completes the proof. \QED

\subsection{Specifying the setup to toric varieties}\label{SS:gs}

Let $M$ be a symplectic toric variety of dimension 
$2n$ and $\De $ denotes its Delzant polytope.
Recall that the  cohomology ring $H^*(M_{T})$ is isomorphic to the face ring
$St(\De )$. In this section,  we give an explicit formula  for fiber integration

$$\pi _*:H^*(M_{T})\cong St(\De)\to H^{*-2n}(B_{T})\cong \B  R[T_1,...,T_n].$$

\NI
In  fact, we give an algorithm which allows to
compute fiber integrals of Chern classes (so also Pontriagin)
of vertical bundle  $TM_{T}$.  The algorithm requires only the data encoded in the
Delzant polytope of $M$.

We  start with few observations. The first is that the Atiyah-Bott
formula become in this case  fairly simple, because the action  has
only isolated fixed points.   Thus  we have

$$\pi  _*(\al ) = \sum _P {i^*_P\al \over E(\nu _P)}\quad ,$$

\noindent
where $i^*_P\al  $ and $E(\nu _P)$ are thought  as ordinary polynomials.

Next we compute  the explicit  form of  the homomorphism
$i^*_P:H^*(M_{T}) \to H^*(B_{T})$ induced by the
inclusion of the fixed point $i_P:{P}\to M$.
Notice that $i^*_P$ is completely
defined by its restriction to $H^2(M_{T})$, 
by multiplicativity.
We describe
it in a matrix  form with  respect  to the following bases:
$\{F_1,F_2,...,F_k\}$ of $H^2(M_{T})$ and $\{T_1,T_2,...,T_n\}$
of $H^2(B_{T})$.  Recall that $F_i$ denotes the  i-th  face of the
Delzant polytope  of  $M$.

Observe that $i^*_P(F) = 0$ for any face $F$ which does not contain the vertex
corresponding to the fixed point $P$. This implies that the column corresponding to
$F$ consists of zeros. Moreover, since  $i_P$ is a section of the universal  bundle
then we have  that $i_P^* \circ \pi ^* = Id$. This two observations
completely determine $i^*_P$. Let's denote the entries of the matrix
representing $i^*_P$ by $a_{ij}^{I_P}$, where $I_P=\{i_1,...,i_n\}$ is such that
$P=F_{i_1}\cap F_{i_2}\cap ...\cap F_{i_n}$.

\BS

\begin{prop}\label{P:matrix}
With the above notation,
the entries of the matrix representing
$i^*_P:H^2(M_{T})\to H^2(B_{T})$ are given by the two following
conditions

\begin{enumerate}
       \item  
       $a^{I_P}_{ij} = 0,\quad   \hbox{if}\quad j\notin I_P$
       \item
       $\sum  _j a_{ij}^{I_P}\cdot u_{jk} = \delta _i^k,$
\end{enumerate}

\noindent
where $u_{jk}$ are the entries of the matrix representing $\pi ^*$.
\end{prop}

\qed

\BS

Finally, recall that
$E(\nu _P) = i^*_Pi_*^P(1) =i^*_P( F_{i_1} F_{i_2}...F_{i_n}) $ \cite{ab},
where  $i_*^P:H^*(B_{T})\to H^{*+2n}(M_{T})$ denotes the push-forward of $i_P$.
Thus this observation together with Proposition \ref{P:matrix} make the
Atiyah-Bott formula explicitly computable.

\section{The computations for symplectic toric varieties}\label{S:main}

\subsection{Rational ruled surfaces}\label{SSS:hs}

The rational cohomology ring of the classifying space of
symplectomorphism group of a rational ruled surface is
known due to \cite{am}. The aim of this subsection
is the ``reality test'', that is to check how big part
of the cohomology ring of $B_{Ham\Mo }$ is generated
by fiber integrals of characteristic classes.

Let $M^k_{\la }$ for $0\leq k\in \B Z$ and
$\la \in \B R$  be a symplectic toric variety, whose Delzant
polytope is a quadrilateral with the following vertices:

$$
\!\!\!\!\!\!\!\!
\left .
\begin{array}{c}
(0,0)\\
(0,1)\\
(1 + \la - {k - 1 \over 2},1)\\
(0,2 + \la + {k - 1 \over 2},0)
\end{array}
\right \} 
\begin{array}{c}
\text{for $k$ odd,}\\
\text{$\la > -1$ and}\\
\text{$2(1 + \la) > k -1$.}
\end{array}
\quad
\text{or}\quad
\left.
\begin{array}{c}
(0,0)\\
(0,1)\\
(1 + \la - {k \over 2},1)\\
(0,1 + \la + {k  \over 2},0)
\end{array}
\right \} 
\begin{array}{c}
\text{for $k$ even,}\\
\text{$\la \geq 0$ and}\\
\text{$2(1 + \la ) > k$.}
\end{array}
$$

\NI

\BS

\NI
\begin{pspicture}(0,0)(3,2)\label{Pic:hs1}
\psgrid [subgriddiv=4,griddots=0,gridlabels=0.15]
\psline[linewidth=1.5pt,fillstyle=solid,fillcolor=lightgray]
       {-}(0,0)(0,1)(1.5,1)(2.5,0)(0,0)
\rput(1.75,-0.5){{\bf Figure 1:}}
\rput(1.5,-0.8){The Delzant polytope of $M^1_{0.5}$}
\end{pspicture}
\qquad \quad
\begin{pspicture}(0,0)(4,2)\label{Pic:hs2}
\psgrid [subgriddiv=4,griddots=0,gridlabels=0.15]
\psline[linewidth=1.5pt,fillstyle=solid,fillcolor=lightgray]
 {-}(0,0)(0,1)(1.5,1)(3.5,0)(0,0)
\rput(1.75,-0.5){{\bf Figure 2:}}
\rput(2,-0.8){ The Delzant polytope of $M^2_{1.5}$}
\end{pspicture}

\BS
\BS
\BS
\BS
\BS
\noindent
These manifolds are symplectic ruled surfaces which means that they
are $S^2$-bundles over $S^2$ and the symplectic form restricts to
symplectic form on each fiber. Symplectic
ruled surfaces (up to rescaling of the symplectic form) are 
symplectomorphic to either $M^0_{\la }$ for some $\la \geq 0$
or $M^1_{\la }$ for some $\la > -1$ \cite{lm1}.
More precisely, we rescale the symplectic forms
such that the area of the fiber is equal to 1.
Then the Delzant polytopes are rectangular trapezoids of  height 1.
Manifolds $M^k_{\la }$ and $M^K_{\La }$ are symplectomorphic
if $k \cong K$ modulo 2 and $\la = \La $.

It follows from the above description that we have defined several 
Hamiltonian torus action on a given symplectic ruled surface
provided that $\la $ is large enough. For example,
manifold $M^2_{1.5}$ (Figure 2) is symplectomorphic to 
$M^0_{1.5}$ and $M^4_{1.5}$. These are three different actions
on one symplectic manifold.

\begin{prop}\label{P:hir}
Let $M^k_{\la }$ for $1 + \la > \lfloor \frac k 2 \rfloor$  
be a symplectic ruled surface
and $H^*(B_T)=\B R[x,y]$. Then the fiber
integrals 
of the characteristic classes are given by the following
formulae.

\begin{enumerate}
\item The monomials in the Chern classes and the coupling class:
\begin{eqnarray*}
\pi_*(c_1^ic_2^j\Om ^l) &=& (x + y)^i (x y)^{j-1} \Om _1^l + \\
& & (x - y)^i (-xy)^{j-1} \Om _2^l + \\
& &((k-1) x - y)^i(-x (k x - y))^{j-1} \Om _3^l+\\
& &((-k-1) x + y)^i (x (k x - y))^{j-1} \Om _4^l,
\end{eqnarray*} 
\NI
where
\begin{eqnarray*} 
\Om _1 &=& {1\over  3 k + 6 \la}\left ( (k^2 + 3 k\la + 3 \la ^2) x +
                   (k + 3 \la ) y \right )\\
\Om _2 &=& {1\over 3 k + 6 \la}\left ( (k^2 + 3 k \la + 3 \la ^2) x -
                   (2 k + 3 \la ) y\right )\\
\Om _3 &=& {1\over 3 k + 6 \la}\left (
                   (k^2 - 3 \la ^2) x -
                   (2 k + 3 \la ) y \right )\\
\Om _4 &=& {1\over 3 k + 6 \la}\left ( 
                   -( 2 k^2 + 6 k \la +  3\la ^2 )x +
                   ( k + 3\la ) y \right ).
\end{eqnarray*}
\item
The fiber integrals of Pontriagin classes are zero.
\item
The fiber integrals of monomials in the Pontriagin and Euler classes:
\begin{eqnarray*}
\pi _*(p_1^i eu^j) \!\!\!\!&=&\!\!\!\!
  (x^2 + y^2)^i \left ((xy)^{j-1} + (-xy)^{j-1}\right ) +\\
&&\!\!\!\!\!\!\!\!\!\!\!\!\!\!\!\!\!\!\!\!\!\! 
((1 + k^2) x^2 - 2kxy + y^2)^i \left 
((x (kx - y))^{j-1} + (-x (kx - y))^{j-1}
\right ).
\end{eqnarray*}
\end{enumerate}
\end{prop}

\NI
{\bf Proof:} This is a direct computation as described in
Section \ref{SS:gs}.\QED

\BS

\begin{cor}\label{C:hir}
Let $\Mo $ be a symplectic ruled surface. 
Then the following statements hold:
\begin{enumerate}
\item
If $\Mo $ is different from
$M^0_1$ and $M^1_{\la }$ for $\la \in (-1,0]$, then
the rational cohomology ring of the classifying
space of the group of symplectomorphisms is generated
by fiber integrals of the monomial in the Chern classes and
the coupling class. 
\item
$H^{4i}(B_{Dif\!f(M)})$ is nontrivial for $i=1,2,3...$
\end{enumerate}
\end{cor}

\NI
{\bf Proof:} 
{\bf (1)} It is known due to Abreu and McDuff \cite{am}, that

$$
H^*(B_{Ham\Mo })= \B Q[A,X,Y]\slash \sim ,
$$

\NI
where $deg(A)=2$, $deg(X)=deg(Y)=4$. The relation is given by a polynomial
depending on $\la $ and the diffeomorphism class of $M$.
Using this result it suffices to check that the fiber integration
is a surjection on the second and fourth cohomology. This can
be done by a direct computation. Indeed,
the above formula for the fiber integrals yields
that $\pi_*^T(c_1^3)\neq 0$ and
$\dim span\{F((\pi_*(c_1^3))^2,F(\pi _*(c_1\Om ^3)),F(\pi _*(c_1^2\Om ^2))\}
=3$, where $F$ denotes the detection function (Definition \ref{D:d}).

{\bf (2)} The nontrivial elements are given by the fiber
integrals of monomials in the Euler and Pontriagin classes.
\QED

\begin{rem}
{\em \hfill

\NI
{\bf 1.} The reason that the first part of the theorem does
not hold in full generality is that the excluded manifolds
don't admit more than one different actions of torus.

\NI
{\bf 2.} Notice that
if $k\neq k^{\prime }$ then the torus actions associated to
$k$ and $k^{\prime }$ are not conjugated in the group of all
diffeomorphisms of $M$. 
If the actions were  conjugate then they would induce
equal maps between $H^*(B_{Dif\! f(M)})$ and $H^*(B_T)$.
This follows from the fact that any conjugation induces an
identity map on the cohomology of the classifying spaces
(see \cite{se}).
Hence
the fiber integrals of the powers of the Euler class would be
equal for different actions. This is a contradiction with the
formula in Proposition \ref{P:hir} (3).

In fact, this can be easily seen directly by
looking at the representations on the tangent spaces
of the fixed points as was pointed  to us by P.Seidel.
This remark also applies to subsequent examples.
}
\end{rem}

\subsection{$\cp ^2\sharp 2 \overline {\cp }^2$}\label{E:5gon}

Our next example consists of the actions on the symplectic manifold
diffeomorphic to
$\cp ^2\sharp 2 \overline {\cp }^2$.
The cohomology ring of the classifying space of
symplectomorphism group is not known in this case.
Our computations provide estimates from below for
its Betti numbers.

Let $\nu \in \B R$ and $k\in \B Z$ be such that
$\nu > 0$ and $0\leq k < \lceil {\nu  \over 2}\rceil $.
Let $M_{\nu ,k}$
be a family of symplectic toric varieties
whose Delzant polytopes are presented on the Figure 2.

\BS
\BS

\centerline{\,\,\,\,\,
\begin{pspicture}(0,0)(8,3)\label{Pic:5gon}
\psset{unit=1.4cm}
\psline[linewidth=1.5pt,fillstyle=solid,fillcolor=lightgray]
       {-}(0,0)(0,1)(1,2)(2,2)(8,0)(0,0)
\rput(-0.25,0.5){\bf $F_1$}
\rput(0.25,1.75){$F_2$}
\rput(1.5,2.25){$F_3$}
\rput(5,1.25){$F_4$}
\rput(5,-0.25){$F_5$}
\rput(-0.1,-0.1){\small 0}
\rput(-0.1,1){\small 1}
\rput(-0.1,2){\small 2}
\rput(1,-0.2){\small 1}
\rput(2,0){\small $\bullet $}
\rput(8,0){\small $\bullet $}
\rput(4,0){\small $\bullet $}
\rput(2,-0.2){\small $\nu + 1 - 2k$}
\rput(8,-0.2){\small $\nu  +1 + 2(k+1)$}
\rput(4,-0.2){\small $\nu + 1 $}
\psgrid [subgriddiv=4,griddots=0,gridlabels=0.0]
\end{pspicture}
}
\BS
\BS
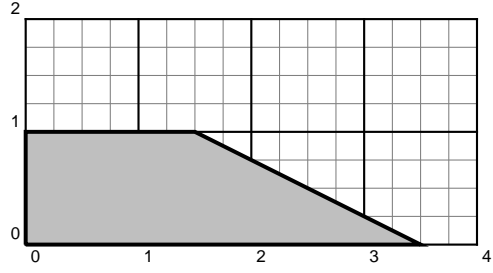
\centerline {\bf Figure 2:}
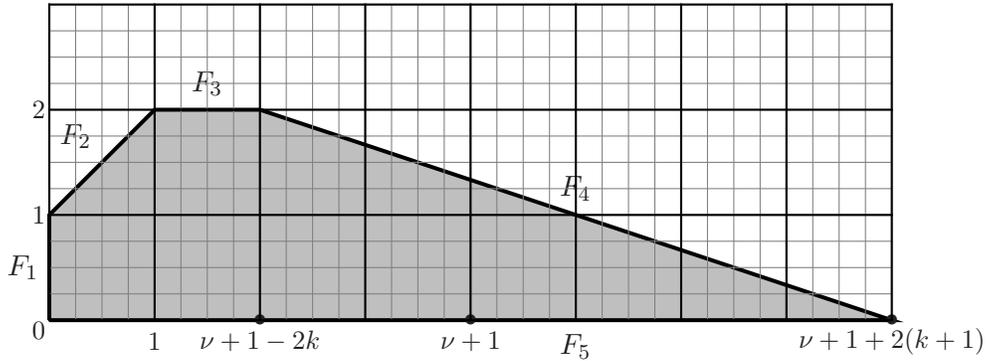
\centerline {The Delzant polytope for a torus action on 
$\cp ^2\sharp 2 \overline \cp ^2$}

\BS
\BS
\NI
Clearly, $M_{\nu ,k}$ are symplectomorphic for fixed $\nu $ and
$0\leq k < \lceil {\nu  \over 2}\rceil $. Indeed, $M_{\nu ,k}$ 
are the symplectic
blow-ups of ruled surfaces $M^{2k+1}_{\nu }$ of the previous example.
Regarding as symplectic manifold they will be denoted by
$M_{\nu }$.

According to  Proposition \ref{P:chern}  we get that
the Chern classes of the vertical bundle correspond to 
edges and vertices and are of the form
$$c_1 = F_1 + ...+ F_5\quad \hbox{and}\quad 
c_2 = F_1F_2 + F_2F_3+ ... + F_5F_1.$$

\noindent
The direct computation gives the following form of
the coupling class

\begin{eqnarray*}
\Om =  {1\over3(7 + 4\nu )}(\! \! \! \! \! \! \! \!  \! \! \! \! &&
          (15 - 8k + 12\nu ) F1 + \\
       & &(25 + 8 k + 8k^2 + 24\nu + 6\nu ^2) F2 + \\
       & &(31 + 16 k + 8k^2 + 24\nu + 6\nu ^2) F3  + \\
       & &(27 + 8k + 12\nu ) F4 + \\
       & &(23 + 12 k + 8 k^2 + 21\nu + 6\nu ^2) F5 ).
\end{eqnarray*}

\begin{prop}\label{P:2cp}
Let $M_{\nu ,k}$ be as above and $H^*(B_T)=\B R[x,y]$. 
Then the fiber
integrals 
of the characteristic classes are given by the following
formulae.

\begin{enumerate}
\item Monomials in the Chern classes and the coupling class:
\begin{eqnarray*}
\pi _*(c_1^ic_2^j\Om ^l) &=& (x + y)^i (x y)^{j-1} \Om_1^l+\\
& & x^i (-y (x + y))^{j-1} \Om_2^l +\\
& & (-y)^i ( -x (x + y))^{j-1} \Om_3^l +\\
& & (2 k  x - y)^i (-x ((2k+1) x - y))^{j-1}\Om_4 ^l+\\
& & (-(2k +2) x - y)^i ( -x (-(2k +1) x + y)^{j-1} \Om_5^l,
\end{eqnarray*}

\NI
where
\begin{eqnarray*}
 \!\!\!\!\!\!\!\!\!\!\!\!\!\!
 \Om _1 \!\!&=& \!\!{1\over 3 (7 + 4 \nu )}\left (
                 (25 + 8k + 8k^2 + 24\nu + 6\nu ^2) x +
                 (15 - 8 k + 12 \nu ) y\right )\\
\!\!\! \!\!\!\!\!\!\!\!\!\!\!
\Om _2  \!\!&=& \!\!{1\over 3 (7 + 4 \nu )} \left (
                 (25 + 8 k + 8 k^2 + 24\nu + 6\nu ^2) x +
                 (- 6 - 8 k) y\right )\\
\!\!\!\!\!\!\!\!\!\!\!\!\!\!
\Om _3  \!\!&=& \!\!{1\over 3 (7 + 4 \nu )}\left (
                 (4 + 8k + 8 k^2 + 12\nu + 6\nu ^2) x +
                 (-27 - 8 k - 12 \nu ) y\right )\\
 \!\!\!\!\!\!\!\!\!\!\!\!\!\!
\Om _4  \!\!&=& \!\!{1\over 3 (7 + 4 \nu )}\left (
                 (4 + 50k + 8k^2 - 9\nu + 24k\nu - 6\nu ^2)x +
                 (-27 - 8 k - 12 \nu ) y\right )\\
\!\!\! \!\!\!\!\!\!\!\!\!\!\!
\Om _5  \!\!&=& \!\!{1\over 3 (7 + 4 \nu )}\left (
                 (-38 - 34 k + 8k^2 - 33\nu - 24k\nu - 6\nu ^2) x +
                 (15 - 8 k + 12 \nu ) y\right ).
\end{eqnarray*}
\item Monomials in the Pontriagin classes:
\begin{eqnarray*}
\pi _*(p_1^ip_2^j) ={1\over xy}&& \!\!\!\!\!\!\!\!\!\!\!
                    \left [(x^2 y^2)^j (x^2 + y^2)^i\right ] -\\ 
           {1\over x y (x + y)}&&\!\!\!\!\!\!\!\!\!\!\!
                \left [y(x^2 (x + y)^2)^j (2 x^2 + 2 x y + y^2)^i + \right .\\
      &&\!\!\!\!\!\!\!\!\!\! \left . 
                 x(y^2 (x + y)^2)^j (x^2 + 2 x y + 2 y^2)^i\right ], 
\end{eqnarray*}
\item Powers of the Euler class:
\begin{eqnarray*}
\pi _*(eu^i) &=& (-x ((2k+1) x - y))^{i-1} + (x ((2k+1) x - y))^{i-1} +\\
              & &(x y)^{i-1} + (-x (x + y))^{i-1} + (-y (x + y))^{i-1}
\end{eqnarray*}

\item Monomials in the Euler and Pontriagin classes:
\begin{eqnarray*}
\pi _*(p_1^i eu^j)\!\!\!\!\!\!\!\!\!\!\!&&= (xy)^{j-1} (x^2 + y^2)^i +\\
\!\!\!\!\!\!\!\!\!\!\!\!&&(-x(x+y))^{j-1} (2x^2 + 2xy + y^2)^i +\\
\!\!\!\!\!\!\!\!\!\!\!\!&&(-y(x+y))^{j-1} (x^2 + 2y(x + y))^i +\\
\!\!\!\!\!\!\!\!\!\!\!\!&&(x(x + 2kx - y))^{j-1} (-2x(x + 2kx - y) + (-2(1+k)x + y)^2)^i +\\
\!\!\!\!\!\!\!\!\!\!\!\!&&(-x(x + 2kx -y))^{j-1} ((2 + 4k + 4k^2)x^2 + y^2 - 2x (y + 2ky))^i.
\end{eqnarray*}

\end{enumerate}\qed
\end{prop}

\begin{cor}\label{C:2cp}
Let $\Mo = M_{\nu }$ be as above.
The following statements are true:
  \begin{enumerate}

\item
If $\nu > 2$, then  $\dim H^2(B_{Ham \Mo }) \geq 4$,

\item
If $0<\nu \leq 2$, then  $\dim H^2(B_{Ham \Mo }) \geq 2$,
      
\item 
$\dim H^4(B_{Diff(M)}) \geq 2$.
Moreover
$H^{4k}(B_{Diff(M)}) $ is nontrivial for $k=1,2,...$.
  \end{enumerate}
\end{cor}

\NI
{\bf Proof:} The first two statements follows from the direct
application of the detection function. 
The assumption on $\nu $ in the first one
ensures that there are at least
two different actions, which contributes to the detection function.
Notice that the first estimate is the best we can get. Indeed, we integrate the
following classes: $\Om ^3,\, \Om^2 c_1,\, \Om c_1^2,\, c_1^3,\,
c_2 \Om,\, c_2c_1 $. $\pi_*(\Om^3)=0$ from the definition and
$\pi _*(c_2c_1) = 0$ according to the $T$-strict multiplicativity
of the Todd genus.

To get the last  statement we integrate $p_1^2$ and $p_1eu$.
This fiber integrals are linearly independent in 
$H^{4}(B_{Diff(M)}) $.
Also 
$\pi _*(p_1^ieu)\neq 0$ which implies that
$H^{4k}(B_{Diff(M)}) $ is nontrivial for $k=1,2,...$.

\qed

\subsection{Dimension 6: projectivizations of complex bundles}

Let $L_k\to \cp^1$ be a line bundle with $c_1(L_k) = k$.
Let $1\leq \mu \in \Bbb R$ be greater than $k$ and $l$.
Consider a  symplectic toric variety
$M_{\mu ,k, l}:=\Bbb P(L_k\oplus L_l\oplus L_0)$,
which is the projectivization of complex bundle
and whose Delzant polytope looks as in the following figure.

\psset{unit=0.75cm}
\centerline{
\begin{pspicture}(0,0)(16,8)\label{Pic:genhir}
\psline[linewidth=1.5pt,fillstyle=solid,fillcolor=lightgray]
       {-}(1,3)(4,1)(10,1)(15,3)(4.5,6)(1,6)(1,3)
\psline[linewidth=1.5pt]{-}(1,6)(4,1)
\psline[linewidth=0.5pt,linestyle=dashed]{-}(10,1)(6.5,3)
\psline[linewidth=1.5pt]{-}(4.5,6)(10,1)
\psline[linewidth=1.5pt]{-}(1,3)(15,3)
\psline[linewidth=0.5pt,linestyle=dashed]{-}(4.5,6)(4.5,3)
\psline[linewidth=1pt]{->}(0.25,3.5)(5.5,0)
\psline[linewidth=1pt]{->}(0.5,3)(16,3)
\psline[linewidth=1pt]{->}(1,2.5)(1,7)
\rput(15,3.4){$\mu $}
\rput(6.5,3.33){$\mu - l$}
\rput(4.5,3.33){$\mu - k$}
\rput(3.7,0.7){1}
\rput(.7,6){1}
\rput(.7,2.7){0}
\end{pspicture}
}
\BS
\centerline {\bf Figure 3}
\centerline {The Delzant polytope for $M_{\mu ,k,l}$}

\BS

\begin{lemma}\label{L:genhir}
Two symplectic toric varieties $M_{\mu ,k,l}$
and 
$M_{\mu ^{\prime },k^{\prime },l^{\prime }}$ 
are symplectomorphic
provided that 

\begin{enumerate}
\item
$k+l \cong k^{\prime } + l^{\prime }\, (\text{mod 3})$ and

\item
$3\mu  - (k+l)= 3\mu ^{\prime } - (k^{\prime } + l^{\prime })$.
\end{enumerate}
\end{lemma}

\pf
First observe that any of the above manifolds is the total space
of a symplectic bundle $\cp ^2\to M \to S^2$, whose structure group
is $PU(3)$. Since $\pi _1(PU(3)) \cong \B Z_3$, then it is clear
that $M_{\mu ,k,l}$ and
$M_{\mu ^{\prime },k^{\prime },l^{\prime }}$  are isomorphic
(as symplectic bundles) if
$k+l \cong k^{\prime } + l^{\prime }\, (\text{mod 3})$.
Moreover, the symplectic structures on these manifolds are
of the form $\Om + K\pi^*(\om _{S^2})$, where
$\Om $ denotes the coupling form and
$\pi:M\to S^2$ is the projection. One can see easily that
$K = \mu - \frac 13(k + l)$ for $M_{\mu ,k,l}$.

Let $M_0$ and $M_1$ be two such isomorphic bundles and 
$H:S^1\times [0,1]\to PU(3)$ be a chosen homotopy
between loops defining $M_0$ and $M_1$. By the usual Thurston
argument, there exist a number $K\in \B R$ such that
$\om _t:=\Om _t + K\pi_t^*(\om _{S^2})$ is a symplectic form
on each $M_t$. 
Notice that $\om _t$ is an isotopy of symplectic form
which does not change the cohomology class, hence it
follows from Moser's argument that $(M_0,\om _0 )$ and
$(M_1,\om _1)$ are symplectomorphic. 

\qed

\begin{rem}{\em
A detailed discussion of the coupling parameter $K$ as well
as exact computations for certain bundles was done by
Polterovich in \cite{po1}.
}
\end{rem}

\BS

Thus the symplectomorphism type depends only on $\mu $ and
on $\la := (k+l)3\Bbb Z  \in \Bbb Z\slash 3\Bbb Z$. 
We denote it by $M_{\mu ,\la }$.
The explicit formulae of fiber integrals are quite complicated, so
we do not present them. As in the previous examples, they allow to estimate
the dimension of cohomology groups.

\begin{prop}\label{P:genhir}
Let $\Mo = M_{\mu ,\la} $ be as above. Then
\begin{enumerate}
\item
$\dim H^2(B_{Ham \Mo }) \geq 1$,

\item
$\dim H^4(B_{Ham \Mo }) \geq 8$, provided that $\Mo $ admits
at least two different Hamiltonian actions of $T^2$.

\item
$H^{4 k}(B_{Diff (M)})$ is nontrivial for $k=1,2,...$. 

\end{enumerate}
\end{prop}

\qed


\section{Restrictions comming from multiplicativity of genera}\label{S:Todd}

When computing fiber integrals one observes many more linear
dependencies between them than we a priori expected. This is
due to multiplicativity properties of certain genera.
Recall that a {\bf genus} is a ring homomorphism

$$
K : \Omega \otimes \B Q \to R,
$$
where $\Omega $ is a cobordism ring and $R$ is an integral domain over
$\B Q$. Since we are working in the symplectic category then
in the sequel $\Omega $ will denote the complex cobordism ring.
Every genus is defined for a stable complex manifold $M$
by a multiplicative sequence $K:=\{K_r \}$, $K_r \in R[x_1,...,x_r]$, as follows.
Given a complex vector bundle $E$ of rank $n$

$$
K(E):= K(c_1(E),...,c_n(E)) := 1 + K_1(c_1(E)) + K_2(c_1(E),c_2(E)) + ...
$$
and

$$
K(M) :=  \left <K(TM);[M]\right > \in R.
$$

\subsection{G-strict multiplicativity}\label{SS:gsm}

Let $G$ be a group.
Multiplicative sequence $K$ is said to be $G$-strictly multiplicative
or $G$-sm for short if 

$$
\pi _*(K(TM_G)) \in H^0(B_G)
$$
for every manifold $M$ on which $G$ acts. Here 
$\pi :M_G \to B_G$ is the universal fibration associated to the
action and $TM_G \to M_G$ is the complex vector bundle tangent
to fibers. In other words, $G$-strict multiplicativity
says that any fiber integral $\pi _*(K_r(c_1(TM_G),...,c_r(TM_G)))=0$
for $r\neq \dim M$. Hence for $r > \dim M$ we obtain relations
in $H^*(B_G)$ between fiber integrals of characteristic classes.

The name multiplicativity is justified by the following fact.

\begin{prop}\label{P:multi}
Let $M\to P \to B$ be a $G$-bundle over a compact stably complex base.
Then $P$ admits a stable complex structure. If $K$ is
a $G$-strict multiplicative sequence then

$$K(P) = K(M) K(B).$$
\end{prop}

\pf
Since $G$-strict multiplicativity is a property of the universal
fibration then it holds for every $G$-fibration. Hence we have that

$$
\pi _*(K(Vert)) \in H^0(B),
$$
where $Vert \to P$ is the bundle tangent to the fibers of
$\pi :P\to B$. The statement follows from the following computation.

\begin{eqnarray*}
K(P)&=& \left <K(TP);[P]\right > = \left <K(Vert\oplus \pi^*TB);[P]\right >\\
&=&\left <K(Vert) K(\pi^*TB);[P]\right >\\
&=&\left <\pi_*[K(Vert)\cup  \pi^*K(TB)];[B]\right >\\
&=&\left <\pi_*K(Vert)\cup  K(TB);[B]\right >\\
&=&\left <K(TM);[M]\right >\left <K(TB);[B]\right >= K(M)K(B)
\end{eqnarray*}

\qed

Notice that if a multiplicative sequence $K$ is $S^1$-sm
then it is also
$G$-sm for any compact Lie group $G$ \cite{och}. Indeed, $H^*(B_G)$
is a subalgebra in $H^*(B_T)$ where $T\subset G$ is a maximal
torus. For any polynomial $p\in \B R[t_1,...,t_n] = H^*(B_T)=\B R[t]$
of degree $k$ (i.e. $p\in H^{2k}(B_T)$)
there exists a homomorphism $f:S^1\to T$ such that
$H^{2k}(B_f)(p) \neq 0$. To see this, notice that any such map is determined
by $n$ integers, say $[k_1,...,k_n]$ and
$H^{2k}(B_f)(p)= p(k_1,...,k_n)t^k$.
Thus if $K$ were not $G$-sm for some
manifold $M$ then it would be not $S^1$-sm.

\subsection{$G$-strict multiplicativity of $\chi_y$-genus}\label{SS:chiy}

Recall that $\chi _y:\Omega \otimes \B Q \to \B Q[y]$ is a genus
whose value on  a K\" ahler manifold is given by

$$
\chi _y(M):= \sum _{p,q} (-1)^q h^{p,q}y^p,
$$
where $h^{p,q}$ are the Hodge numbers \cite{hbj}.

\begin{thm}\label{T:chiy}
$\chi _y$-genus is $G$-strictly multiplicative for every
compact Lie group $G$.
\end{thm}

\pf
According to the observation in the previous subsection, it 
is sufficient to prove strict multiplicativity
for $S^1$. This is equivalent to the fact that $\chi_y$-genus
is multiplicative for any fibration
$M\to P\to \cp ^k$, that is

$$
\chi_y(P) = \chi_y(M)\chi_y(\cp ^k).
$$

We compute $\chi_y(P)$ with the use of the localization formula
for circle action. The formula due to Kosniowski and Lusztig allows
to compute the genus in terms of the genera of the fixed points:

$$
\chi_y(P) = \sum _{i=1}^m (-y)^{s_i}\chi_y(F_i).
$$
Here $F_i$ is the component of the fixed point set and 
$s_i := \dim V_{>0}$ the dimension of the subspace of the
tangent space to any point of $F_i$ on which the circle
acts positively ($z\cdot v = z^k v$, where $k>0$).

We start with defining a circle action on $P$. Let
$\al :S^1 \to Aut(M)$ be a given action on $M$ preserving
the stable complex structure. 
Moreover let $\be :S^1 \to Aut(\cp ^k)$ be an action given by
$\be (t)[z_0:...:z_k] = [t^{i_0}z_0:...:t^{i_k}z_k]$ with
isolated fixed points. We use the same notation for
the action lifted to $S^{2k+1}$.
The bundle $P$ is associated
to the principal bundle $S^{2k+1}\to \cp^k$, i.e.

$$
P:= S^{2k+1}\times _{S^1} M,
$$
where $(z,x)\simeq (z_1,x_1)$ iff $z_1 = s z$ and $x_1 = \al (s^{-1})(x)$,
for $s\in S^1$.

We define an action $\phi :S^1 \to Aut(P)$ by

$$
\phi (t)[z,x] = [\be (t)(z),\al (t)(x)].
$$
Clearly, this action is well defined since $S^1$ is commutative.
Further, fixed points lie in the fibers over the fixed points
of $\be $. When restricted to the fixed fiber, the fixed points of
$\phi $ are those of $\al $. More precisely, let $F_{j}$ denote
a path component of the fixed point set of $\al $. Denote by $F_{ij}$
the $j$-th component of the fixed point det
of $\phi $ lying in the fiber over the $i$-th fixed point of $\be$.
Clearly, $F_{ij} \cong F_j$.  

We need to figure out the dimension $s_{ij}$ of the subspace of the tangent
space to any fixed point $x\in F_{ij}$ on which the circle acts positively.
Notice that the infinitesimal action in the direction normal to
the fiber is the same as infinitesimal action induced by $\be $ on
the tangent space to the fixed point in $\cp ^k$.
Thus $s_{ij}=s_i + s_j$, where $s_i$ ($s_j$ respectively) denote the appropriate dimension
with respect to the action $\be $ ($\al $ respectively).

Now, plugging the above observations into the localization formula
we get the statement as follows.

\begin{eqnarray*}
\chi _y(P) 
&=& \sum _{i=1}^{k+1}\sum _j (-y)^{s_{ij}}\chi_y(F_{ij})\\
&=& \sum _{i=1}^{k+1}\sum _j (-y)^{s_i + s_j}\chi_y(F_{ij})\\
&=& \sum _j\left (\sum _{i=1}^{k+1}(-y)^{s_{i}}(-y)^{s_j}\right)\chi_y(F_j)\\
&=& \left(\sum _j (-y)^{s_j}\chi_y(F_j)\right) \left(\sum _{i=1}^{k+1}(-y)^{s_i}\right)\\
&=& \chi_y(M)\chi_y(\cp ^k).
\end{eqnarray*}

\qed

\begin{cor}\label{C:todd}
The Todd genus and the signature  is $G$-sm for any compact Lie group.
\end{cor}

\pf
The Todd genus is equal to $\chi _y$-genus for $y=0$ and the signature
for $y=1$.

\qed

It would be very interesting to know to what extent strict
multiplicativity is true. For example, 
the signature is $Ham\Mo $-sm,
since the group $Ham\Mo $
of Hamiltonian symplectomorphisms
is connected.
These considerations motivate for the following

\BS
\NI
{\bf Question:} Is the Todd (or $\chi _y-$) genus
$Ham\Mo $-strictly multiplicative
for any compact symplectic manifold?

\BS

\begin{ex}[The Todd genus is not $Symp-$strictly multiplicative.]\hfill
{\em
Let $\Si_h \to M \to \Si_g$ be a symplectic surface bundle 
over surface with
nonzero signature, $\si (M)\neq 0$ \cite{a}. 

\begin{eqnarray*}
3 \si (M)
\!\!\!\! &=&\!\!\!\! 
\left <p_1(TM),[M]\right > = \left <c_1(TM)^2 - 2 c_2(TM),[M]\right > \\
\!\!\!\! &=&\!\!\!\! 
\left <c_1(Vert \oplus \pi^*T\Si_g)^2,[M]\right > - 
2\left  <c_2(Vert \oplus \pi^*T\Si_g),[M]\right >  \\
\!\!\!\! &=&\!\!\!\! 
\left <c_1(Vert)^2 + 2c_1(Vert)\pi^*c_1(T\Si_g),[M]\right > -
2\left  <c_1(Vert)\pi^*c_1(T\Si_g),[M]\right > \\
\!\!\!\! &=&\!\!\!\! 
\left <(c_1(Vert)^2),[M]\right > \\
\!\!\!\! &=&\!\!\!\! 
\left <\pi_*(c_1(Vert)^2),[\Si_g]\right > \neq 0
\end{eqnarray*}

\NI
Since the Todd polynomial $T_2(c_1,c_2) = \frac 1 {12} (c_1^2 + c_2)$,
then the Todd genus is not strictly multiplicative in this case.
In this example signature measures the defect of strict multiplicativity.
}
\end{ex}

\NI
{\bf Authors addresses:}\\

\NI
$
\begin{array}{lll}
\text {Tadeusz Januszkiewicz}       &\quad  &\text {Jarek K\c edra}\\
\text {Mathematical Institute UWr}  &\quad  &\text {Institute of Mathematics US}\\
\text {pl. Grunwaldzki 2/4}        &\quad  &\text {Wielkopolska 15}\\
\text {50-384 Wroc\l aw}           &\quad &\text {70-451 Szczecin}\\
\text {Poland}                     &\quad  &\text {Poland }
\end{array}
$

\end{document}